\documentclass[12pt]{article}

\usepackage{amsmath}
\usepackage{amssymb}
\usepackage{amsthm}

\newtheorem{theorem}{Theorem}

\newtheorem{proposition}{Proposition}

\title{Powered numbers in short intervals II}
\author{Tsz Ho Chan}
\date{}

\begin{document}
\maketitle

\begin{abstract}
In this article, we derive better results concerning powered numbers in short intervals, both unconditionally and conditionally on the $abc$-conjecture. We make use of sieve method, a polynomial identity, and a recent breakthrough result on density of sets with no $k$-term arithmetic progression. In the process, we study integers over short intervals that have with a big smooth divisor.
\end{abstract}

\section{Introduction and main results}

A number $n$ is {\it squarefull} or {\it powerful} if its prime factorization $n = p_1^{a_1} p_2^{a_2} \cdots p_r^{a_r}$ satisfies $a_i \ge 2$ for all $1 \le i \le r$. Similarly, $n$ is $k$-{\it full} if $a_i \ge k$ for all $1 \le i \le r$. In contrast, $n$ is {\it squarefree} if $a_i = 1$ for all $i$. For example, $72 = 2^3 \cdot 3^2$ is squarefull, $648 = 2^3 \cdot 3^4$ is $3$-full, and $30 = 2 \cdot 3 \cdot 5$ is squarefree. For any integer $n$, we define its squarefree part $q(n)$ to be
\[
q(n) := \prod_{i \text{ with } a_i = 1} p_i = \prod_{p \mid n, \, p^2 \nmid n} p
\]
and any integer can be factored uniquely into the product of its squarefree part and its powerful part (for example, $q(72) = 1$, $q(360) = 5$, and $360 = 5 \cdot 2^3 \cdot 3^2 = 5 \cdot 72$.) In \cite{C:2023}, the author studied powerful numbers in short intervals and obtained
\[
\mathop{\sum_{x < n \le x + y}}_{n \text{ is powerful}} 1 \ll \frac{y}{\log (y +2)},
\]
and
\begin{equation} \label{smooth}
\mathop{\mathop{\sum_{x < n \le x + y}}_{n \text{ is powerful}}}_{p^{+}(n) \le \sqrt{y} } 1 \ll y^{11/12}
\end{equation}
for $1 \le y \le x$ where $p^{+}(n)$ denotes the largest prime factor of $n$. Similarly, we let $p^{-}(n)$ denote the smallest prime factor of $n$.

\bigskip

Generalizing and smoothing $k$-full numbers, Mazur \cite{M:2000} proposed the study of {\it powered numbers} in connection to the $abc$-conjecture: For any $\epsilon > 0$, there exists a constant $C_\epsilon > 0$ such that, for any integers $a, b, c$ with $a + b = c$ and $\text{gcd}(a,b) = 1$, the bound
\[
\max\{ |a|, |b|, |c| \} \le C_\epsilon \kappa(a b c)^{1 + \epsilon}
\]
holds where $\kappa(n) := \prod_{p \mid n} p$ is the {\it squarefree kernel} of $n$. 

Let $k \ge 1$ be a real number (not necessarily an integer). A positive integer $n$ is {\it $k$-powered} if
\[
\kappa(n) \le n^{1/k} \; \; \text{ or } \; \; \frac{\log n}{\log \kappa(n)} \ge k.
\]
For example, $648$ is a $\frac{\log 648}{\log 6} = 3.613147...$-powered number. One can easily see that any $k$-full number is also $k$-powered (so the latter is a more general notion). Conditional on the $abc$-conjecture, the author \cite{C:2024} proved that, with $k > 3/2$,
\begin{equation} \label{kpow}
\mathop{\sum_{x < n \le x + y}}_{n \text{ is $k$-powered}} 1 \ll_k \frac{y}{\exp( (c_k \log y)^{0.09})}
\end{equation}
for $1 \le y \le x$; and, for integers with smooth squarefree part,
\begin{equation} \label{verysmooth}
\mathop{\sum_{x < n \le x + y}}_{p^+(q(n)) \le \log (y + 1) \log \log (y + 2)} 1 \ll \frac{y \log \log (y +2)}{\log (y + 2)}.
\end{equation}

In this note, we obtain some new results in connection with (\ref{smooth}), (\ref{kpow}) and (\ref{verysmooth}).
\begin{theorem} \label{thm1}
Let $k > 1$ and $0 \le \delta < 1$ be any real numbers. For $1 \le y \le x$,
\[
\mathop{\mathop{\sum_{x < n \le x + y}}_{n \text{ is $k$-powered}}}_{p^{+}(n / \kappa(n)) \le y^{1 - \delta}} 1 \; \ll_k \; y^{\frac{3k + 1}{4k}} + \frac{y^{1 - \delta}}{\log (y + 1)}.
\]
\end{theorem}
The condition $p^{+}(\frac{n}{\kappa(n)}) \le y^{1 - \delta}$ means that the powerful part of $n$ is $y^{1 - \delta}$-smooth. Hence, Theorem \ref{thm1} generalizes (\ref{smooth}) to powered numbers, and gives the better exponent $7/8$ (when $\delta = 1/2$ and $k = 2$) instead of $11/12$ in (\ref{smooth}).

\begin{theorem} \label{thm2}
Let $k > 1$ be any real number. For $2 \le w \le y \le x$ and $\log w = o(\log y)$,
\[
\mathop{\mathop{\sum_{x < n \le x + y}}_{n \text{ is $k$-powered}}}_{p^+(q(n)) \le w} 1 \; \ll \; \frac{y \log w}{\log y} + y \exp \Bigl(- \frac{\log y}{3 \log w} \Bigr) + y^{11/12}.
\]
\end{theorem}
If one focuses on powered numbers only, then Theorem \ref{thm2} improves on (\ref{verysmooth}) by allowing for a much wider range of smoothness on the squarefree part of $n$.

\begin{theorem} \label{thm3}
Assume the $abc$-conjecture. For $k > 5/4$, there exist positive constants $C_k$ and  $c$ such that
\[
\mathop{\sum_{x < n \le x + y}}_{n \text{ is $k$-powered}} 1 \ll_k \frac{y}{\exp( C_k (\log \log y)^{c})}.
\]
for $1 \le y \le x$.
\end{theorem}
In comparison with (\ref{kpow}), Theorem \ref{thm3} extends the range of $k$ to $k > 5/4$ in exchange for an inferior upper bound. It would be interesting to see if other polynomial identities or ideas could be used to extend the range all the way to $k > 1$.

\bigskip

The paper is organized as follows. First, we will use Shiu's method \cite{S:1980} to prove Theorem \ref{thm1}. Then, we will prove an upper bound on the number of integers with a relatively big smooth divisor over short intervals (which may be of independent interest), and apply it to deduce Theorem \ref{thm2}. Finally, we will use a recent breakthrough result by Leng, Sah and Sawhney \cite{LSS:2024} on density of sets with no $k$-term arithmetic progression to show Theorem \ref{thm3}.

\bigskip

Throughout the paper, $p$, $q$, $p_i$ and $q_j$ stand for prime numbers. The symbols $f(x) = O(g(x))$, $f(x) \ll g(x)$, and $g(x) \gg f(x)$ are equivalent to $|f(x)| \leq C g(x)$ for some constant $C > 0$. Also, $f(x) = O_{\lambda} (g(x))$, $f(x) \ll_{\lambda} g(x)$, and $g(x) \gg_\lambda f(x)$ mean that the implicit constant may depend on the parameter $\lambda$. The symbol $f(x) = o(g(x))$ means that $\lim_{x \rightarrow \infty} \frac{f(x)}{g(x)} = 0$. The function $\exp(x)$ is defined to be $\exp(x) := e^{x}$. For any set $\mathcal{S}$, $| \mathcal{S}|$ denotes the number of elements in $\mathcal{S}$. Finally, a number is {\it $y$-smooth} if all of its prime divisors are less than or equal to $y$.

\section{Smooth powerful part: Theorem \ref{thm1}}

\begin{proof}[Proof of Theorem \ref{thm1}]
We may assume that $y > 1$ for the theorem is clearly true for $y = 1$. We assume $z \le y^{(k - 1) / k}$. Any $k$-powered number $n$ in $[x, x + y]$ can be factored uniquely as $n = a b$ with $a$ squarefree and $b$ squarefull with $\gcd(a, b) = 1$. Suppose we have the following prime factorizations:
\[
a = \underbrace{p_1 \cdots p_i}_{a_1} \underbrace{p_{i+1} \cdots p_r}_{a_2}\; \; \text{ with } \; \; p_1 < p_2 < \cdots < p_r,
\]
and
\[
b = \underbrace{q_1^{e_1} \cdots q_j^{e_j}}_{b_1} \underbrace{q_{j+1}^{e_{j+1}} \cdots q_s^{e_s}}_{b_2} \; \; \text{ with } \; \; q_1 < q_2 < \cdots < q_s \le y
\]
where $i$ and $j$ are the greatest indices such that $p_1 \cdots p_i \le y / z$ and $q_1^{e_1} \cdots q_j^{e_j} \le z$. Hence, $a_1 \le y / z < a_1 p_{i+1}$ and $b_1 \le z < b_1 q_{j+1}^{e_{j+1}}$. Note that $i$ or $j$ may be $0$ (giving empty products for $a_1$ or $b_1$) if $p_1 > y / z$ or $q_1^{e_1} > z$. Also, $i$ and $j$ may equal to $r$ and $s$ respectively.

\bigskip

Case 1: $b_1 > z^{1/2}$. The number of such $k$-powered numbers is bounded by
\begin{equation} \label{c1}
\mathop{\sum\nolimits'}_{z^{1/2} < b \le z} \; \mathop{\sum_{x < n \le x + y}}_{b \mid n} 1 \le  \mathop{\sum\nolimits'}_{z^{1/2} < b \le z} \frac{y}{b} \ll \frac{y}{ z^{1/4}}
\end{equation}
where $\mathop{\sum\nolimits'}$ denotes a sum over squarefull numbers only.

\bigskip

Now, we claim that if $b_1 \le z^{1/2}$, then $b_2$ must be greater than $1$. Suppose the contrary that $n = a_1 a_2 b_1$. If $\kappa(a_1)^k \kappa(b_1)^k > a_1 b_1$, then, as $n$ is $k$-powered and $\kappa(a b) = \kappa(a) \kappa(b)$ when $\gcd(a, b) = 1$, we have
\[
n = a_1 a_2 b_1 \ge \kappa(a_1)^k \kappa(a_2)^k \kappa(b_1)^k = \kappa(a_1)^k \kappa(b_1)^k a_2^k > a_1 b_1 a_2^k
\]
which is a contradiction. Hence, $\kappa(a_1)^k \kappa(b_1)^k \le a_1 b_1$ (i.e., $a_1 b_1$ itself is a $k$-powered number) and 
\[
n = a_1 a_2 b_1 \ge a_1^k a_2^k \kappa(b_1)^k \; \; \text{ or } \; \; a_2^{k - 1} \le \frac{b_1}{a_1^{k - 1} \kappa(b_1)^k}.
\]
Since $n > x$, this implies
\[
\Bigl( \frac{x}{a_1 b_1} \Bigr)^{k-1} < \frac{b_1}{a_1^{k - 1} \kappa(b_1)^k} \; \; \text{ or } \; \; y \le x < \Bigl( \frac{b_1}{\kappa(b_1)} \Bigr)^{k / (k-1)} \le z^{\frac{k}{2(k-1)}}
\]
as $b_1 / \kappa(b_1) \le b_1 \le z^{1/2}$. This contradicts our choice $z \le y^{(k-1) / k}$.

\bigskip

Case 2: $b_1 \le z^{1/2}$ and $p_{-}(b_2) \le z^{1/2}$. Then $p_{j+1} \le z^{1/2}$ and $p_{j+1}^{e_{j+1}} > z^{1/2}$ which implies $p_{j+1}^{-e_{j+1}} \le \min(z^{-1/2}, p_{j+1}^{-2})$ as $e_{j+1} \ge 2$. Hence, the sum
\[
\sum_{p_{j+1} \le z^{1/2}} \frac{1}{p_{j+1}^{e_{j+1}}} \le \sum_{p_{j+1} \le z^{1/4}} z^{-1/2} + \sum_{z^{1/4} < p_{j+1} \le z^{1/2}} \frac{1}{p_{j+1}^2} \ll \frac{1}{z^{1/4}}.
\]
Therefore, by replacing $p_{j+1}^{a_{j+1}}$ with a generic $p^a$, the number of $k$-powered numbers in this case is bounded by 
\begin{equation} \label{c2}
\sum_{p \le z^{1/2}} \; \mathop{\sum_{x < n \le x + y}}_{p^a \mid n} 1 \le \sum_{p \le z^{1/2}} \Bigl( \frac{y}{p^a} + 1 \Bigr) \ll \frac{y}{z^{1/4}} + z^{1/2}.
\end{equation}

Case 3: $b_1 \le z^{1/2}$ and $z^{1/2} < p_{-}(b_2) \le y^{1 - \delta}$. In this case, there is some prime $z^{1/2} < p \le y^{1 - \delta}$ such that $p^2$ divides $n$. The number of such $k$-powered numbers is bounded by
\begin{equation} \label{c3}
\sum_{z^{1/2} < p \le y^{1 - \delta}} \; \mathop{\sum_{x < n \le x + y}}_{p^2 \mid n} 1 \le \sum_{z^{1/2} < p \le y^{1 - \delta}} \Bigl( \frac{y}{p^2} + 1 \Bigr) \ll \frac{y}{z^{1/4}} + \frac{y^{1 - \delta}}{\log y}.
\end{equation}
Combining (\ref{c1}), (\ref{c2}) and (\ref{c3}), we have Theorem \ref{thm1} by taking $z = y^{(k-1) / k}$.
\end{proof}

\section{Smooth squarefree part: Theorem \ref{thm2}}

First, let us recall a standard sieve upper bound over short intervals: For any $x$ and any $y \ge 2$,
\begin{equation} \label{sievebd}
\mathop{\sum_{x < n \le x + y}}_{\gcd(n, \, \prod_{p \le \sqrt{y}} p ) = 1} 1 \le \frac{2y}{\log y} \Bigl(1 + O \Bigl( \frac{1}{\log y} \Bigr) \Bigr).
\end{equation}
For a proof, see Theorem 3.3 of \cite{MV:2007} for example. Next, we deduce an upper bound on the number of integers in short intervals that have with a relatively big smooth divisor. This result is an adaptation of Lemma 2.2 of \cite{BKS:2021} to short intervals.
\begin{proposition} \label{P}
Let $0 < \alpha \le 1/2$ and $\log w = o(\log y)$ with $y$ sufficiently large. Consider
\[
\mathcal{S}_\alpha := \{ x < n \le x + y : \text{ there exists $d \mid n$ such that } p^+(d) \le w \text{ and } d > y^\alpha \}.
\]
Then
\[
| \mathcal{S}_\alpha | \ll y \Bigl( \exp \Bigl(- \alpha \frac{\log y}{\log w} \Bigr) + y^{-\alpha / 3} \Bigr).
\]
Note: When $\alpha > 1/2$, one can simply observe that $\mathcal{S}_\alpha \subset \mathcal{S}_{1/2}$ and get the upper bound $|\mathcal{S}_\alpha| \ll y \exp( -\frac{ \log y}{2 \log w} ) + y^{5/6}$.
\end{proposition}

\begin{proof}
Since $\log w = o (\log y)$, we may assume that $w \le y^{\alpha / 4}$. For $n \in \mathcal{S}$, let $p_1 \le p_2 \le \cdots \le p_k$ be the sequence of prime factors of $n$ of size $\le w$ listed in increasing order and according to their multiplicity. From the definition of $\mathcal{S}$, we have $p_1 \cdots p_k > y^\alpha$. Let $j$ be the smallest index such that $p_1 \cdots p_j > y^\alpha$. We must have $j \ge 5$ as $p_i \le y^{\alpha / 4}$. Now, set
\[
a = p_1 \cdots p_{j-2}, \; \; \; p = p_{j-1}, \; \; \text{ and } \; \; b = \frac{n}{a p}.
\]
Hence, $a > y^\alpha /(p_{j-1} p_j) \ge y^{\alpha / 2}$, $a p \le y^\alpha$, and $p^+(a) \le p \le p_{-}(b)$. Therefore,
\begin{equation} \label{Scount}
| \mathcal{S}_\alpha | \le \sum_{p \le w} \, \mathop{\sum_{y^{\alpha / 2} < a \le y^\alpha / p}}_{p^+(a) \le p} \, \mathop{\sum_{x/(ap) < b \le (x + y) /(a p)}}_{p_{-}(b) \ge p} 1 \ll \sum_{p \le w} \, \mathop{\sum_{y^{\alpha / 2} < a \le y^\alpha / p}}_{p^+(a) \le p} \frac{y}{a p \log p}
\end{equation}
by sieve bound (\ref{sievebd}) since $a p \le y^{\alpha} \le \sqrt{y}$ yielding $y / (a p) \ge \sqrt{y} \ge p$. If we let $\epsilon_p = \min \{ \frac{2}{3}, \frac{2}{\log p} \}$, then Rankin's trick yields
\[
\frac{| \mathcal{S}_\alpha |}{y} \ll \sum_{p \le w} \, \mathop{\sum_{y^{\alpha / 2} < a \le y^\alpha / p}}_{p^+(a) \le p} \frac{(a / y^{\alpha / 2})^{\epsilon_p}}{a p \log p} \le \sum_{p \le w} \frac{y^{- \alpha \epsilon_p / 2}}{p \log p} \sum_{p^+(a) \le p} \frac{1}{a^{1 - \epsilon_p}} \ll \sum_{p \le w} \frac{y^{- \alpha \epsilon_p / 2}}{p \log p} \prod_{q \le p} \Bigl(1 - \frac{1}{q^{1 - \epsilon_p}} \Bigr)^{-1}
\]
with $q$ denoting a prime number. Since $q^{\epsilon_p} = 1 + O(\log q / \log p)$ for $q \le p$, Mertens' estimates imply the inside product is $\ll \log p$. Hence,
\begin{align*}
\frac{| \mathcal{S}_\alpha |}{y} &\ll y^{- \alpha / 3} + \sum_{100 < p \le w} \frac{e^{- \alpha \log y / \log p}}{p} \le y^{- \alpha / 3} + \sum_{j \ge 1} \sum_{w^{1 / (j+1)} < p \le w^{1 / j}} \frac{e^{- j \alpha \log y / \log w}}{p} \\
&\ll y^{-\alpha / 3} + \sum_{j \ge 1} e^{- j \alpha \log y / \log w} \ll y^{-\alpha / 3} + e^{- \alpha \log y / \log w}
\end{align*}
by Mertens' estimates again. This completes the proof of Proposition \ref{P}.
\end{proof}

\begin{proof}[Proof of Theorem \ref{thm2}]

This proof is similar to that of Theorem \ref{thm1}. Hence, we will highlight the necessary modifications at times. First, we may assume that $y > c_0$ is sufficient large as the theorem is clearly true for $y \le c_0$ by picking a sufficiently large implicit constant. We leave $z$ as a free parameter and factor $n = a b \in (x, x +y]$, $a$ and $b$ in the same way as the previous section. 

\bigskip

First, suppose $4 \log \log y \le \log w = o( \log y )$. We set $z = w$. Then, inequality (\ref{c1}) from case 1 yields an upper bound $O(y / z^{1/4}) = O(y / \log y)$. When $b_1 \le z^{1/2}$, the same argument in Theorem \ref{thm1} shows that $b_2 > 1$. Then, inequality (\ref{c2}) in case 2 gives an upper bound $O(y / z^{1/4}) = O(y / \log y)$. The main difference is that $p_{-}(b_2) \le y$ may not hold for case 3. 

\bigskip

Case 3(a): $b_1 \le \sqrt{z}$ and $\sqrt{z} < p_{-}(b_2) \le \sqrt{y}$. Then, inequality (\ref{c3}) yields an upper bound $O(y / z^{1/4}) = O(y / \log y)$.

\bigskip

Case 3(b): $b_1 \le \sqrt{z}$ and $p_{-}(b_2) > \sqrt{y}$. Recall $a_1 \le y / z < a_1 p_{i+1}$. Then $a_1 b_1 \le y / \sqrt{z} \le y$. Since the squarefree part $a$ is $w$-smooth, we must have $a_1 > y / z w$ unless $a = a_1 \le y / z w$. If $y / z w < a_1 \le y / z$, we can apply Proposition \ref{P} with $\alpha = 1/2$. This implies that the number of such $k$-powered numbers is bounded by $O(y \exp(- \frac{\log y}{2 \log w} ) + y^{5/6})$. It remains to deal with the case $a \le y / z w$. 

\bigskip

Subcase 1: $a > \sqrt{y / z}$. By Proposition \ref{P} with $\alpha = 1/3$, the number of such $k$-powered numbers is at most $O(y \exp( - \frac{\log y}{3 \log w}) + y^{11/12})$.

\bigskip

Subcase 2: $a \le \sqrt{y / z}$. Then the number of such $k$-powered numbers is at most
\[
\sum_{a \le \sqrt{y/z}} \; \sum_{b_1 \le z^{1/2}} \mathop{\sum_{\frac{x}{a b_1} < b_2 \le \frac{x + y}{a b_1}}}_{p_{-}(b_2) > \sqrt{y}} 1 \ll \sum_{a \le \sqrt{y/z}} \; \sum_{b_1 \le z^{1/2}} \frac{y}{a b_1 \log y} \ll \frac{y}{\log y} \prod_{p \le w} \Bigl(1 + \frac{1}{p} \Bigr) \ll \frac{y \log w}{\log y}
\]
by (\ref{sievebd}) (as $\frac{y}{a b_1} > \sqrt{y}$) and Mertens' estimate. Combining all of the above bounds, we have Theorem \ref{thm2} when $4 \log \log y \le \log w = o(\log y)$.

\bigskip

Now, we assume that $\log w \le 4 \log \log y$. We set $z = \exp(\sqrt{\log y})$. Then case 1, case 2 and case 3(a) above yield the bound $O(y / z^{1/4}) = O(y / \log y)$. For case 3(b), when $y /z^2 < a_1 \le y / z$, we have the bound $O(y^{5/6} + y \exp(- \frac{\log y}{8 \log \log y})) = O( y / \log y )$ by Proposition \ref{P} with $\alpha = 1/2$. 

\bigskip

It remains to deal with the case $a \le y / z^2$. Similar to subcase 1 above, when $a > \sqrt{y / z}$, the number of such $k$-powered numbers is at most $O(y^{11/12} + y \exp(- \frac{\log y}{12 \log \log y})) = O( y / \log y )$ by Proposition \ref{P} with $\alpha = 1/3$.

\bigskip

Similar to subcase 2 above, when $a \le \sqrt{y / z}$, the number of such $k$-powered numbers is at most
\[
\sum_{a \le \sqrt{y/z}} \; \sum_{b_1 \le z^{1/2}} \mathop{\sum_{\frac{x}{a b_1} < b_2 \le \frac{x + y}{a b_1}}}_{p^+(b_2) > \sqrt{y}} 1 \ll \sum_{a \le \sqrt{y/z}} \; \sum_{b_1 \le z^{1/2}} \frac{y}{a b_1 \log y} \ll \frac{y}{\log y} \prod_{p \le w} \Bigl(1 + \frac{1}{p} \Bigr) \ll \frac{y \log w}{\log y}.
\]
Combining all of the above bounds, we have Theorem \ref{thm2} when $\log w \le 4 \log \log y$.
\end{proof}

\section{Conditional result: Theorem \ref{thm3}}

Let $r_k(N)$ denote the size of the largest subset of $\{1, 2, \ldots, N \}$ with no non-trivial $k$-term arithmetic progressions. Here, non-trivial means that the $k$-terms are not the same. Now, let us recall Gowers' quantitative breakthrough result on Szemer\'{e}di's theorem \cite{G:2001}:
\[
r_k(N) \le N (\log \log N)^{-c_k} \; \; \text{ with } \; \; c_k = 2^{-2^{k+9}}.
\]
For $k = 3$, there have been recent exciting activities \cite{BS:2020}, \cite{KM:2023}, \cite{BS:2023a}, \cite{BS:2023b} in getting down to
\[
r_3(N) \ll N \exp( - c (\log N)^{1/9}).
\]
For $k = 4$, Green and Tao \cite{GT:2009}, \cite{GT:2017} proved
\[
r_4(N) \ll N (\log N)^{-c}
\]
for some $c > 0$. For $k \ge 5$, Leng, Sah, and Sawhney \cite{LSS:2024} very recently established the existence of some constant $c_k > 0$ such that
\begin{equation} \label{lss}
r_k(N) \ll N \exp (- (\log \log N)^{c_k} ).
\end{equation}

\begin{proof}[Proof of Theorem \ref{thm3}]
When $k > 5/4$, we set $\epsilon := k / 5 - 1 / 4 > 0$. Firstly, we restrict to $y \le x^{\frac{4k - 5}{8k}}$. We claim that there is no non-trivial $7$-term arithmetic progression of $k$-powered numbers in the interval $(x, x + y]$ under the $abc$-conjecture. Suppose the contrary that
\[
n - 2d, \; n - d, \; n, \; n + d, \; n + 2d, \; n + 3d, \; n + 4d
\]
are seven $k$-powered numbers in $(x, x + y]$ with some $d \ge 1$. We apply the polynomial identity
\[
(n + 2d)^3 (n - 2d) + 16 d^3 (n + d) = n^3 (n + 4d). 
\]
Let $t := \gcd(n, d)$, $n := t n'$, and $d := t d'$ for some positive integers $n'$ and $d'$ with $\gcd(n', d') = 1$. Then
\[
(n' + 2d')^3 (n' - 2d') + 16 d'^3 (n' + d') = n'^3 (n' + 4d').
\]
Say $D := \gcd( n'^3 (n' + 4 d'), 16 d'^3 (n' + d'))$. Suppose some prime power $p^l \mid D$. Note that $p \nmid d'$ for otherwise $p \mid d'$ and $p \mid D \mid n'^3 (n' + 4 d')$ would imply $p \mid n'$ by Euclid's lemma. This contradicts $\gcd(d', n') = 1$. Since $p^l \mid D \mid 16 d'^3 (n' + d')$, we must have $p \mid 16$ or $p \mid n' + d'$ by Euclid's lemma.

\bigskip

Case 1: $p > 3$. We have $p \mid n' + d'$. Since $p \mid D \mid n'^3 (n' + 4 d')$, we have $p \mid n'$ or $p \mid n' + 4 d'$ by Euclid's lemma. The former situation implies $p \mid (n' + d') - n' = d'$ which contradicts $\gcd(d', n') = 1$. The latter situation implies $p \mid (n' + 4 d') - (n' + d') = 3 d'$. This forces $p \mid d'$ and $p \mid (n' + d') - d' = n'$ which contradicts $\gcd(d', n') = 1$ again. Therefore, this case cannot happen.

\bigskip

Case 2: $p = 3$. Then we have $3^l \mid n' + d'$ and $3^l \mid n'^3 (n' + 4d')$. By Euclid's lemma, we have $3 \mid n'$ or $3 \mid n' + 4 d'$. The former situation cannot happen for otherwise it would imply $3 \mid (n' + d') - n' = d'$ and $3 \mid n'$. Hence, we must have $3^l \mid n' + 4 d'$ and $3^l \mid (n' + 4 d') - (n' + d') = 3 d'$. If $l \ge 2$, then $3 \mid d'$ and $3 \mid (n' + d') - d' = n'$, a contradiction. Therefore, $l = 1$.

\bigskip

Case 3: $p = 2$. Since $2^l \mid D \mid n'^3 (n' + 4 d')$, we have $2 \mid n'$ and $2 \nmid n' + d'$ as $2 \nmid d'$. Thus, we must have $2^l \mid 16$.

\bigskip

Summarizing the above cases, we have $D = 2^{e_2} 3^{e_3}$ for some integers $0 \le e_2 \le 4$ and $0 \le e_3 \le 1$. In particular, $1 \le D \le 48$. Now, we apply the $abc$-conjecture to
\[
\frac{(n' + 2d')^3 (n' - 2d')}{D} + \frac{16 d'^3 (n' + d')}{D} = \frac{n'^3 (n' + 4d')}{D}
\]
and get
\begin{align*}
\frac{x^4}{48 t^4} &\ll_\epsilon \kappa\Bigl(\frac{(n' + 2d')^3 (n' - 2d')}{D} \Bigr)^{1 + \epsilon} \kappa\Bigl(\frac{16 d'^3 (n' + d')}{D} \Bigr)^{1 + \epsilon} \kappa\Bigl(\frac{n'^3 (n' + 4d')}{D} \Bigr)^{1 + \epsilon} \\
&\le \kappa ( (n' + 2d')^3 (n' - 2d') )^{1 + \epsilon} \kappa (16 d'^3 (n' + d') )^{1 + \epsilon} \kappa( n'^3 (n' + 4d'))^{1 + \epsilon} \\
&\le \bigl( \kappa(n - 2d) \kappa(n) \kappa(n + d) \kappa(n + 2d) \kappa(n + 4d) \bigr)^{1 + \epsilon} \Bigl( \frac{16 d}{t} \Bigr)^{1 + \epsilon} \ll x^{\frac{5 (1 + \epsilon)}{k}} \cdot \frac{d^{1 + \epsilon}}{t^{1 + \epsilon}}
\end{align*}
by $\kappa(a) \le \kappa(a b) \le \kappa(a) \kappa(b)$, $\kappa(a^3) = \kappa(a)$, and $\kappa(a) \le a$. As $t \le d \le y$, the above inequality implies
\[
x^{4 - \frac{5 (1 + \epsilon)}{k}} \ll_\epsilon d^{1+\epsilon} t^{3 - \epsilon} \le y^4 \; \; \text{ or } \; \; y \gg_\epsilon x^{1 - \frac{5 (1 + \epsilon)}{4 k}}.
\]
From the definition of $\epsilon$, we have $y \ge D_k x^{(12 k - 15)/(16 k)}$ for some constant $D_k > 0$ which contradicts our assumption $y \le x^{(4k - 5)/(8k)}$ when $x > D_k^{16k/(4k - 5)}$. This completes the proof of no non-trivial $7$-term arithmetic progression of $k$-powered numbers in $(x, x + y]$ when $x > D_k^{16k/(4k - 5)}$. Since arithmetic progressions are invariant under translation, we may shift all the $k$-powered numbers in $(x, x+y]$ to numbers in $(0, y]$ and there is no $7$-term arithmetic progression among them. By (\ref{lss}), the number of such $k$-powered numbers is bounded by $C \cdot y \exp( - (\log \log y)^{c_7})$ for some $C \ge 1$ and $c_7 > 0$. When $y \le x^{(4k - 5)/(8k)}$ and $x \le D_k^{16k/(4k - 5)}$, we have $y \le D_k^2$. The number of $k$-powered numbers in $(x, x + y]$ is trivially bounded by $y \le C \cdot \exp((\log \log D_k^2)^{c_7}) \cdot y \exp( - (\log \log y)^{c_7})$. Hence, in any case, we have
\[
\mathop{\sum_{x < n \le x + y}}_{n \text{ is $k$-powered}} 1 \le C \cdot \exp((\log \log D_k^2)^{c_7}) \cdot \frac{y}{\exp( (\log \log y)^{c_7})}
\]
when $y \le x^{(4k - 5) / (8k)}$.

\bigskip

When $x^{(4k - 5)/(8k)} < y \le x$, we can cover the interval $(x, x + y]$ by a disjoint union of at most $2 y / x^{(4k - 5)/(8k)}$ subintervals, each with length $x^{(4k - 5)/(8k)}$. Over each subinterval, we have at most
\[
C \cdot \exp((\log \log D_k^2)^{c_7}) \cdot \frac{x^{(4k - 5)/(8k)}}{ \exp ( (\log \log x^{(4k - 5)/(8k)})^{c_7} ) } 
\]
$k$-powered numbers. Adding all these together, the total number of $k$-powered numbers in $(x, x + y]$ is bounded by 
\[
2 C \cdot \exp((\log \log D_k^2)^{c_7}) \cdot \frac{y}{ \exp ( C_k (\log \log y)^{c_7} )}
\]
for some constant $C_k > 0$.
\end{proof}

\noindent {\bf Acknowledgement.} The author would like to thank the anonymous referee for very helpful suggestions.

%---------------------------------------------------------------------------------------------------------

\end{document}